\documentclass[12pt]{amsart}
 
\input{xy} 
\xyoption{all}

\newtheorem{theorem}{Theorem} [section]
\newtheorem{prop}[theorem]{Proposition}

\newtheorem*{question}{Question}

\numberwithin{equation}{section} 
\usepackage{graphicx}

\newcommand\C{{\mathbb C}}

\renewcommand\P{{\mathbb P}}
\newcommand\R{{\mathbb R}}
\newcommand\Z{{\mathbb Z}}
\newcommand\N{ {\mathbb N}}

\newcommand\cM{\mathcal{M}}
\newcommand\cC{\mathcal{C}}
\newcommand{\cT}{{\mathcal T}}

\newcommand{\cG}{\mathcal {G}}

\newcommand\eps{\varepsilon}

\renewcommand\phi{\varphi}

\renewcommand\mod{\operatorname{mod}}  

\newcommand\MP {\mathrm{MP}}

\newcommand\ord {\operatorname{ord}} 
\newcommand\TF {\operatorname{TF}}
\newcommand\SF {\operatorname{SF}}
\newcommand\Spines {\operatorname{Spines}}
\newcommand\Top {\operatorname{Top}}

\begin{document}

\title{Enumerating the basins of infinity of cubic polynomials}

\author{Laura De Marco and Aaron Schiff}

\begin{abstract}  
We study the dynamics of cubic polynomials restricted to their basins of infinity, and we enumerate topological conjugacy classes with given combinatorics.
\end{abstract}

\date{\today}
\thanks{Research of both authors supported by the National Science Foundation.}

\maketitle

\thispagestyle{empty}

\section{Introduction}

Let $f: \C\to\C$ be a polynomial of degree 3 with complex coefficients.  Its {\em basin of infinity} is the open invariant subset 
	$$X(f) = \{z\in\C: f^n(z)\to\infty \mbox{ as } n\to\infty\}.$$
In this article, we examine combinatorial topological-conjugacy invariants of the restricted dynamical system $f : X(f) \to X(f)$ and count the number of possibilities for each invariant.  We apply results of \cite{DP:combinatorics} which show that we can use these invariants to classify topological conjugacy classes of pairs $(f, X(f))$ within the space of cubic polynomials.  Moreover, when $f$ is in the {\em shift locus}, meaning that both of its critical points lie in $X(f)$, these invariants classify conjugacy classes of polynomials $f:\C\to\C$.  

Specifically, we implement an algorithm which counts topological conjugacy classes of cubic polynomials of {\em generic level N}, defined by the condition that 
	$$G_f(c_1)/3^N < G_f(c_2) < G_f(c_1)/3^{N-1};$$
here $\{c_1, c_2\}$ is the set of critical points of the polynomial $f$ and 
	$$G_f(z) = \lim_{n\to\infty} \frac{1}{3^n} \log^+|f^n(z)|$$
is the escape-rate function.  These generic cubic polynomials are precisely the structurally stable maps in the shift locus \cite{McS:QCIII} (see also \cite{DP:heights}). 
We examine the growth of the number of these stable conjugacy classes as $N\to\infty$. 

We begin with an enumeration of the Branner-Hubbard tableaux (or equivalently, the Yoccoz $\tau$-functions) of length $N$, as introduced in \cite{Branner:Hubbard:2}; see Theorem \ref{tau extension}.  Using the combinatorics of tableaux, we provide an algorithm for computing the number of {\em truncated spines} (introduced in \cite{DP:combinatorics}) for each $\tau$-function; see Theorem \ref{spine count}.  Finally, we apply the procedure of \cite{DP:combinatorics} to count the  number of generic topological conjugacy classes associated to each truncated spine.  
The ideas and proofs follow the treatment of cubic polynomials in \cite{Branner:Hubbard:1}, \cite{Branner:Hubbard:2}, \cite{Branner:cubics}, and \cite{Blanchard:Devaney:Keen}.

\subsection{Results of the computation.}
An implementation of the algorithm was written with Java.  
We compiled the output in Table \ref{table} to level $N = 21$, together with run times (Processor: 2.39GHz Intel Core 2 Duo,  Memory: 2 $\times$ 1GB PC2100 DDR 266MHz). 

\begin{table}
\begin{center}  \begin{small}
\begin{tabular}{ | c | c | c | c | c | c | }
	\hline
Level & Tau sequences & Trees  & Truncated spines & Conjugacy classes & Run time \\  \hline
1	& 	1	&	1 	&	1	&	1 	&	0.000 \\
2	&	2	&	2 	& 	2 	&	2  	&	0.078 \\
3 	&	4 	&	4 	&	4 	&	4 	&	0.062 \\
4 	&	8 	&	8 	&	8 	&	8 	&	0.063  \\
5 	&	16 	&	18	&	18 	&	19 	&	0.093 \\ 
6	&	33 	&	42	&	42 	&	46 	&	0.079  \\
7 	&	69 	&	103	&	105 	&	118 &	0.078  \\
8 	&	144 &	260	&	270 	&	318 &	0.093  \\
9 	&	303 &	670	&	718 	&	881 &	0.094 \\
10 	&	641 &	1753	&	1939 	& 2480 	&	0.125 \\
11 	&  1361 	&	4644	&	5312 	& 7084 	&     0.156  \\
12	&  2895	&	12433	&	14719	& 20374	&     0.266   \\
13	&  6174	&	33581	&	41161	& 59061	&     0.547   \\
14	& 13188	&	91399	&  115856	& 172016	&     1.141  \\
15	& 28229	&	250452	&  328098	& 503018	&     2.453  \\
16 	& 60515	&	690429	&  933719	& 1475478 	&     5.515  \\
17	& 129940 & 1913501	& 2668241	& 4338715	&   12.500  \\
18	& 279415	&	& 7652212	& 12785056	&	27.109  \\
19	& 601742	&	& 22013683	& 37739184	&	72.579  \\
20	& 1297671	&	& 63497798 	& 111562926 &	163.422  \\
21	& 2802318	&	& 183589726 &  330215133  & 383.640  \\  \hline
\end{tabular} \bigskip
\caption{Enumeration of conjugacy classes to generic level $N=21$, with run times measured in seconds.  The tree numbers were computed in [DM].}
\label{table} 
\end{small} \end{center}  
\end{table}

\begin{table}
\begin{center}  \begin{small}
\begin{tabular}{ | c | c | c | c | c | }
	\hline
Levels 17 / 16 & Levels 18 / 17	& 	Levels 19 / 18	&	Levels 20 / 19 &   Levels 21 / 20  \\
2.941 &  2.947	&	2.952	&	2.956 	& 	2.960 	 \\  \hline
\end{tabular} \bigskip
\caption{Computing the growth:  ratios of numbers of conjugacy classes at consecutive levels.}
\label{table2} 
\end{small} \end{center}  
\end{table}

\subsection{The tree of cubic polynomials.}
\begin{figure} 
\includegraphics[width=3.0in]{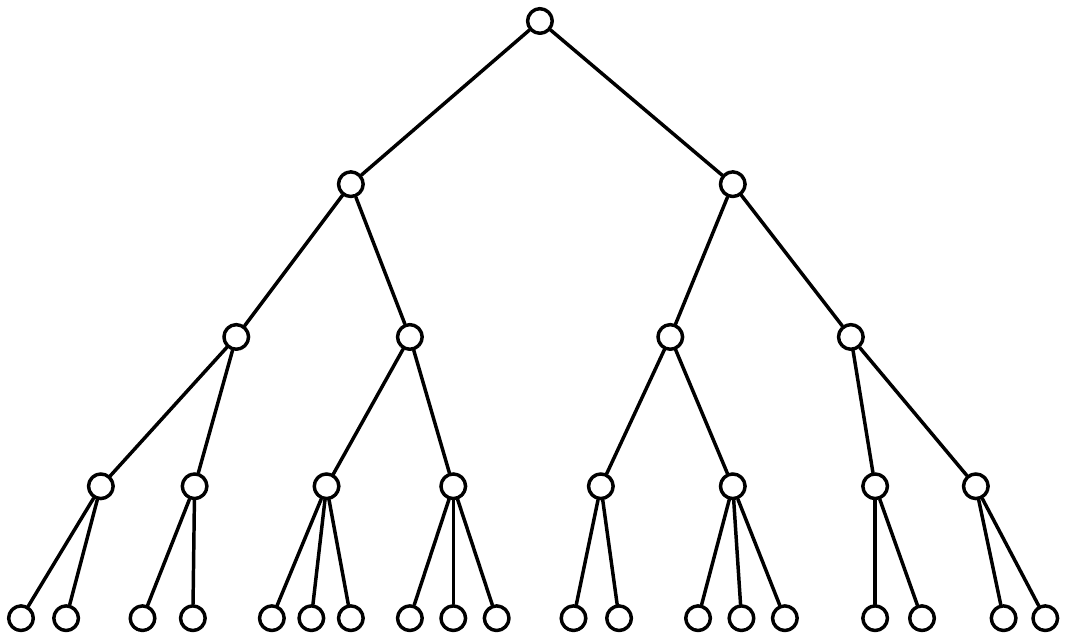}
\caption{The tree $\cT_{conj}$ of conjugacy classes to generic level $N=5$.} \label{tree figure}
\end{figure}

We can define a tree $\cT_{conj}$ of conjugacy classes of cubic polynomials in the shift locus as follows.  For each generic level $N\geq 1$, let $V(N)$ be a set of vertices consisting of one vertex for each topological conjugacy class.  We connect a vertex $v$ in $V(N+1)$ to a vertex $w$ in $V(N)$ by an edge if for each representive $f$ of $v$ and $g$ of $w$, there exists an $\eps>0$ so that the restrictions $f| \{G_f > (G_f(c_1)/3^{N-1})  - \eps\}$ and $g| \{G_g > (G_g(c_1)/3^{N-1}) - \eps\}$ are topologically conjugate.  In Figure \ref{tree figure}, we have drawn this tree $\cT_{conj}$ to level $N=5$.  The growth rate of the number of conjugacy classes as $N\to\infty$ corresponds to a computation of the ``entropy" of this tree. 

The tree $\cT_{conj}$ can also be constructed in the following way.  
Let $\MP_3$ denote the space of conformal conjugacy classes of cubic polynomials; it is a two-dimensional complex orbifold with underlying manifold isomorphic to $\C^2$.  Every cubic polynomial is conformally conjugate to one of the form
	$$f(z) = z^3 + az + b,$$
which can be represented in $\MP_3$ by $(a, b^2)\in\C^2$.  

The critical escape-rate map
	$$\cG: \MP_3\to\R^2$$
is defined by $\cG(f) = (G_f(c_1), G_f(c_2))$ where the critical points $\{c_1, c_2\}$ of $f$ are labeled so that $G_f(c_1)\geq G_f(c_2)$; it is continuous and proper \cite{Branner:Hubbard:1}.  The fiber of $\cG$ over the origin in $\R^2$ is the connectedness locus $\cC_3$, the set of polynomials with connected Julia set.  If we restrict $\cG$ to its complement $\MP_3\setminus\cC_3$, there is an induced projectivization:
	$$\bar{\cG}: \MP_3\setminus\cC_3 \to [0,1]$$
defined by $f\mapsto G_f(c_2)/G_f(c_1)$.  The quotient space of $\MP_3\setminus \cC_3$ formed by collapsing connected components of fibers of $\bar{\cG}$ to points is a (completed) tree $\P\cT_3^*$:  over $(0,1]$ it forms a locally finite simplicial tree while $\bar{\cG}^{-1}(0)$ forms its space of ends. 

In \cite{DP:heights}  it is proved that the edges of the tree $\P\cT_3^*$ correspond to generic topological conjugacy classes.  Thus, the combinatorial tree $\cT_{conj}$ is {\em dual} to the shift locus tree in $\P\cT_3^*$; each edge in $\P\cT_3^*$ corresponds to a vertex in $\cT_{conj}$.  

The tree $\P\cT_3^*$ comes equipped with a projection to the {\em space of cubic trees} $\P\cT_3$ introduced in \cite{DM:trees}.  In \cite{DM:trees} the growth of the number of edges in $\P\cT_3$ was studied (see Table \ref{table}, third column), but the value of the entropy was left as an open question.  Furthermore, the tree $\cT_{conj}$ is a quotient of the tree of {\em marker autorphisms} $\cM_3$  introduced in \cite{Blanchard:Devaney:Keen}; the quotient is by the monodromy action from a twisting deformation (see \cite{Branner:cubics}).  The entropy of $\cM_3$ was easily shown to be $\log 3$, and so the entropy of $\cT_{conj}$ (or equivalently, of $\P\cT_3^*$) is no more than $\log 3$.  

\begin{question}
Let $\phi(N)$ denote the number of vertices in $\cT_{conj}$ at level $N$.  Is 
	$$\lim_{N\to\infty} \frac{1}{N} \log \phi(N) = \log 3 \; ? $$ 
\end{question}

\medskip
In Table \ref{table2}, we show the ratios of the number of conjugacy classes in consecutive levels.  As the level increases, the computed ratios increase, conjecturally limiting on 3.


\bigskip

\section{The $\tau$ functions}

\subsection{The $\tau$-function of a polynomial.} \label{tau definition}
Fix a cubic polynomial $f$ with disconnected Julia set, and let $c_1$ and $c_2$ be its critical points, labeled so that $G_f(c_2)\leq G_f(c_1)$.   For each integer $n$ such that $G_f(c_2) < G_f(c_1)/3^{n-1}$, we define the {\em critical puzzle piece} $P_n(f)$ as the connected component of $\{z: G_f(z) < G_f(c_1)/3^{n-1}\}$ containing $c_2$, and set	$$\tau_f(n) = \max\{j< n: f^{n-j}(c_2) \in P_j(f)\}.$$ 

Recall that the {\em tableau} or {\em marked grid} of $f$ is an array $\{M_f(j,k) \in \{0,1\}: j,k\geq 0\}$, defined by the condition 
	$$M_f(j,k) = 1 \iff f^k(c_2) \in P_j(f).$$
We depict a marked grid as a subset of the 4th quadrant of the $\Z^2$-lattice, where $j\geq 0$ represents the distance along the negative $y$-axis and $k\geq 0$ respresents the distance along the positive $x$-axis.  The values of $\tau_f$ can be read directly from the marked grid:  beginning with $M_f(n,0)=1$, $\tau_f(n)$ is the $j$-coordinate at the first non-zero entry when proceeding ``northeast" from $(n,0)$.  In fact, the orbit $\{\tau_f^k(n): k\geq 0\}$ consists of the $j$-coordinates of all non-zero entries along the diagonal $M_f(n-i, i)$.  Thus, the marked grid can be recovered from the $\tau$-function by:
	$$M_f(j,k) = \left\{ \begin{array}{ll}
			1	& 	\mbox{if } j = k = 0 \\
			1	& 	\mbox{if } j = \tau_f^m(j+k) \mbox{ for some } m\geq 0 \\
			0	& 	\mbox{otherwise} \end{array} \right. $$

Branner and Hubbard \cite[Theorem 4.1]{Branner:Hubbard:2} showed that marked grids associated to cubic polynomials are characterized by a simple set of rules.  A marked grid of {\em size} $N$ (which may be infinite) is an array $\{M(j,k)\in\{0,1\}: j,k\geq 0 \mbox{ and } j+k\leq N\}$ which satisfies the following rules:  
\begin{itemize}
\item[(M0)] 	For each $n\leq N$, $M(n,0) = M(0,n) = 1$.
\item[(M1)]	If $M(j,k) = 1$, then $M(l,k) = 1$ for all $l\leq j$.
\item[(M2)]	If $M(j,k) = 1$, then $M(j-i,k+i) = M(j-i,i)$ for all $0 \leq i \leq j$.
\item[(M3)]	If $j+k<N$, $M(j,k) = 1$, $M(j+1, k)=0$, $M(j-i,i)=0$ for $0< i< m$, and $M(j-m+1,m)=1$, then $M(j-m+1, k+m) = 0$.
\item[(M4)]	If $j+k< N$, $M(j,k)=1$, $M(1,j)=0$, $M(j+1,k)=1$, and $M(j-i,k+i)=0$ for all $0<i<j$, then $M(1,j+k)=0$.
\end{itemize}
The rule (M4) was omitted in \cite{Branner:Hubbard:2}, though it is necessary for their proof. It appears as stated here in \cite[Proposition 4.5]{Kiwi:cubics}; an equivalent formulation (in the language of $\tau$-functions) was given in \cite{DM:trees}.

\subsection{Properties of tau-functions}  \label{properties}
Let $\N$ denote the positive integers $\{1, 2, 3, \ldots\}$.
We consider the following five properties of functions $\tau: \N \to \N\cup\{0\}$.
\begin{itemize}
\item[(A)]	$\tau(1) = 0$
\item[(B)]	$\tau(n+1) \leq \tau(n)+1$
\end{itemize}
From (A) and (B), it follows that $\tau(n) < n$ for all $n\in\N$; consequently, there exists a unique integer $\ord(n)$ such that the iterate $\tau^{\ord(n)}(n) = 0$.  
\begin{itemize}
\item[(C)]	If $\tau(n+1) < \tau^k(n) + 1$ for some $0 < k < \ord(n)$, 
			then $\tau(n+1) \leq \tau^{k+1}(n) + 1$.
\item[(D)]	If $\tau(n+1) < \tau^k(n) + 1$ for some $0 < k < \ord(n)$, 
			and if $\tau(\tau^k(n) + 1) = \tau^{k+1}(n) + 1$, 
			then $\tau(n+1) < \tau^{k+1}(n) + 1$.
\item[(E)]	If $\ord(n)>1$ and $\ord(\tau^{\ord(n)-1}(n) + 1)=1$, then $\tau(n+1)\not=0$.
\end{itemize}

\begin{prop} \label{tau properties}
For any positive integer $N$, a function 
	$$\tau: \{1, 2, 3, \ldots, N\} \to \N\cup\{0\}$$
or a function 
	$$\tau: \N \to \N\cup\{0\}$$
is the $\tau$-function of a cubic polynomial if and only if it satisfies properties (A)--(E).
\end{prop}

\noindent
We say the $\tau$-function is {\em admissible} if it satisfies properties (A)--(E).

The proof is by induction on $N$.  It is not hard to see that the $\tau$ function must satisfy these rules, by doing a translation of the tableau rules.  Conversely, any tau function satisfying properties (A)-(E) determines a marked grid satisfying the 4 tableau rules.  Property (E) is another formulation of the ``missing tableau rule" (M4) appearing in [Ki] and [DM].  

\subsection{Algorithm to inductively produce all $\tau$-functions}  \label{tau algorithm}
If a $\tau$-function has domain $\{1, 2, 3, \ldots, N\}$, we say it has {\em length N}.  The {\em markers} of a $\tau$ with length $N$ are the integers 
	$$\{m\in \{1, \ldots, N-1\}: \tau(m+1) < \tau(m)+1\}.$$
Let $k$ be the number of markers which appear in the orbit 
	$$N\mapsto \tau(N) \mapsto \ldots \mapsto \tau^{\ord(N)}(N)=0,$$ 
and label these $k$ markers by $l_1', l_2', \ldots, l_k'$ so that
	$$N = l_0' > l_1' > l_2' > \cdots > l_k' > 0.$$
For each $0\leq i \leq k$, let $l_i = \tau(l_i')$ so that 
	$$\tau(N) = l_0 > l_1 > \cdots > l_k \geq 0.$$

\begin{theorem} \label{tau extension}
Given an admissible $\tau$-function of length $N$, an extension to length $N+1$ is admissible if and only if 
	$$\tau(N+1) = l_i + 1 \mbox{ for some } 0 \leq i \leq k$$
or $\tau(N+1)=0$ if $l_k>0$ or $k=0$.  
\end{theorem}

\proof
The theorem follows from Proposition \ref{tau properties}.  Property (C) implies that $\tau(N+1)$ must be either 0 or of the form $\tau^k(N) +1$.  Property (D) implies that $\tau(N+1)$ must be either 0 or of the form $l_i +1$.  Property (E) implies that $\tau(N+1)\not=0$ if $k>0$ and $l_k=0$.  On the other hand, if $k=0$, both $\tau(N+1) = l_0+1 = \tau(N)+1$ and $\tau(N+1) = 0$ are admissible. 
\qed

\bigskip
\section{The truncated spine} \label{spine}

Suppose $f$ is a cubic polynomial of generic level $N$.  Introduced in \cite{DP:combinatorics}, the {\em truncated spine} of $f$ is a combinatorial object which carries more information than the $\tau$-sequence though it does not generally determine the topological conjugacy class.  (It determines the tree of local models for $f$, studied in \cite{DP:combinatorics}.)  Here, we describe how to inductively construct truncated spines, and we compute the number of extensions to a truncated spine of length $N+1$ from one of length $N$.  We show exactly how many distinct extensions correspond to a choice of $\tau$-function extension.  

\subsection{The truncated spine of a polynomial}
Fix a cubic polynomial $f$ of generic level $N$.  The truncated spine is a sequence of $N$ finite hyperbolic laminations, one for each connected component of the critical level sets of $G_f$ separating the critical value $f(c_1)$ from the critical point $c_2$, together with a labeling by integers $<N$.  

Specifically, beginning with the level $G_f(c_1)$ of fastest-escaping critical point, we identify the level set $\{G_f = G_f(c_1)\}$ with the quotient of a metrized circle.  The level curve is topologically a figure 8, metrized by its external angles, giving it total length $2\pi$.  The critical point $c_1$ lies at the singular point of the figure 8, identifying points at distance $2\pi/3$ along a metrized circle.  Thus the associated hyperbolic lamination in the unit disk consists of a single hyperbolic geodesic joining two boundary points at distance $2\pi/3$.   The lamination is only determined up to rotation.  We mark the complementary component in the disk with boundary length $4\pi/3$ with a $0$ to indicate the component of $\{G_f < G_f(c_1)\}$ containing the second critical point $c_2$.  

For each critical level set $\{G_f = G_f(c_1)/3^n\}$, $0 < n < N$, we consider only the connected component which separates $c_2$ from $c_1$.  External angles can be used to metrize the curve, normalizing the angles so the connected curve has total length $2\pi$ for each $n$.  The curve can thus be represented as a quotient of a metrized circle; the associated hyperbolic lamination consists of a hyperbolic geodesic joining each pair of identified points.    The {\em gaps} of the lamination are the connected components of the complement of the lamination in the disk; each gap corresponds to a connected component of $\{G_f < G_f(c_1)/3^n\}$.  A gap in a lamination is labeled by the integer $k\geq 0$ if the corresponding component contains $f^k(c_2)$.  

\begin{figure} 
\includegraphics[width=2.5in]{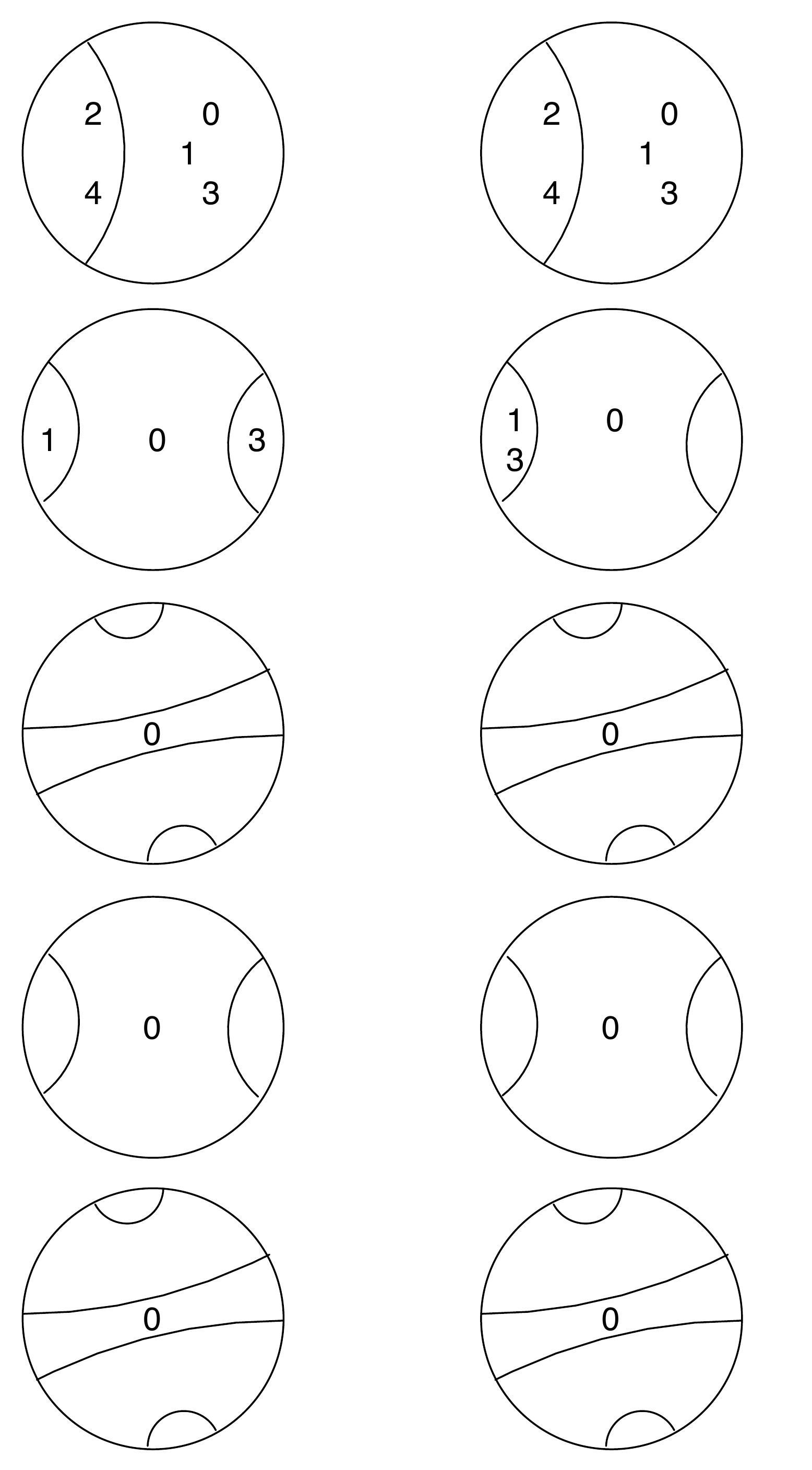}
\caption{Two truncated spines associated to polynomials of generic level $N=5$ with the same $\tau$-sequence $0, 1, 0, 1, 0$.} \label{ex1}
\end{figure}

In Figure \ref{ex1}, we provide examples of truncated spines for two polynomials of generic level $N=5$ with the same $\tau$-sequence.

Two truncated spines of length $N$ are {\em equivalent} if for each $0 \leq n < N$, the labeled laminations at level $n$ are the same, up to a rotation.  In particular, the labeling of the gaps must coincide.  It is shown in \cite{DP:combinatorics} that a truncated spine (together with the heights $(G_f(c_1), G_f(c_2))$ of the critical points) determines the full tree of local models for a polynomial of generic level $N$.  In particular, it carries more information than the {\em tree} of $f$, introduced in \cite{DM:trees}, and its $\tau$-function; it does {\em not}, however, determine the topological conjugacy class of the polynomial.  

The {\em central component} of a lamination is the gap containing the symbol $0$; it corresponds to the component containing the critical point $c_2$.  All non-central gaps are called {\em side components}.  

\subsection{From truncated spine to $\tau$-function}
Fix a truncated spine of length $N$.  Recall that every $\tau$-function satisfies $\tau(1) = 0$.  For $N>1$ and each $0 < n < N$, we can read $\tau(n)$ directly from the truncated spine by: 
	$$\tau(n) = \max\{j: \mbox{lamination } j \mbox{ is labeled by } (n-j)\}.$$
To compute $\tau(N)$, we consider the set
	$$L(N) = \{j : \mbox{the central component in lamination } j \mbox{ is labeled by } (N - j - 1)\}.$$
We then have 
	$$\tau(N) = \left\{ \begin{array}{ll}
		1 + \max\{j: j\in L(N) \} & \mbox{ if } L(N) \not= \emptyset \\
		0	& \mbox{ if } L(N) = \emptyset  \end{array} \right. $$

\subsection{Extending a spine of length $N$}
Fix a truncated spine of length $N$.  It follows directly from the definitions that an extension to length $N+1$ is completely determined by the location of the label $N-\tau(N)$ at level $\tau(N)$.  Any choice of central or side component is admissible:  it determines the local model for an extended tree of local models (from \cite{DP:combinatorics}).   

The lamination at level $N$ is then constructed by taking a degree 2 branched cover of the lamination at level $\tau(N)$ branched over the gap containing the label $N-\tau(N)$.  See \cite{DP:combinatorics} for a general treatment of branched covers of laminations.  The labels are added inductively:  for each iterate $1 < n \leq \ord(N)$, the label $(N-\tau^n(N))$ is placed in the gap at level $\tau^n(N)$ which is the image of the gap containing $\tau^{n-1}(N)$ at level $\tau^{n-1}(N)$.  

\subsection{Computing the number of extended spines for each choice of $\tau(N+1)$}  
Fix a truncated spine with its $\tau$-function of length $N$.  As in \S\ref{tau algorithm}, the {\em markers} of $\tau$ are the integers 
	$$\{m\in \{1, \ldots, N-1\}: \tau(m+1) < \tau(m)+1\}.$$
The {\em marked levels} of $\tau$ are all integers in the forward orbits of the markers:
	$$\{l\geq 0: l = \tau^n(m) \mbox{ for marker } m \mbox{ and } n>0\} \cup \{0\};$$
we say 0 is marked even if there are no markers.    

As before, we let $k$ be the number of markers which appear in the orbit 
	$$N\mapsto \tau(N) \mapsto \ldots \mapsto \tau^{\ord(N)}(N)=0.$$ 
Label these $k$ markers by $l_1', l_2', \ldots, l_k'$ so that
	$$N = l_0' > l_1' > l_2' > \cdots > l_k' > 0.$$
For each $0\leq i \leq k$, let $l_i = \tau(l_i')$ so that 
	$$\tau(N) = l_0 > l_1 > \cdots > l_k \geq 0.$$
For each $0 \leq i < k$, define $n_i$ be the condition that 
	$$\tau^{n_i} (l_i) = l_{i+1}$$
and define $n_k$ so that $\tau^{n_k}(l_k)=0$.  
For $0 < i < j \leq k+1$, we set 
	$$\delta(i,j) = \left\{ \begin{array}{ll} 
		1	&	\mbox{if } \tau(l_i'+1) = l_j+1 \\
		0	&	\mbox{otherwise} \end{array} \right. $$
where by convention we take $l_{k+1} = -1$.   Note that $\tau(l_k'+1) = 0$ for every $\tau$, so $\delta(k, k+1)=1$.  

At level $l_0 = \tau(N)$, there are $2^{\ord(l_0)} = 2^{n_0 + n_1 + \cdots + n_k}$ side components and one central component.  Labeling the central component with the integer $(N-\tau(N))$ uniquely corresponds to the choice of $\tau(N+1) = \tau(N)+1$.  For each $i>0$, the number of side components which correspond to the choice $\tau(N+1) = l_i + 1$ is 
	$$2^{n_0} ( 2^{n_1} ( 2^{n_2} ( \cdots (2^{n_{i-1}} - \delta(i-1,i)) - \cdots ) - \delta(2, i) ) - \delta(1,i) );$$
as above, we take $l_{k+1} = -1$.  It remains to consider how many distinct truncated spines these side components determine.

The {\em symmetry} of $\tau$ is 
	$$s = \min\{n\geq 0: \tau^n(l_0) \mbox{ is a marked level}  \}.$$
Note that $s\leq n_0$.   To each admissible choice for $\tau(N+1)$ (from Theorem \ref{tau extension}) we define the $(N+1)$-th {\em spine factor} of $\tau$.  If $\tau(N+1) = l_i + 1$, then 
$$ \SF(\tau, N+1) = \left\{ \begin{array}{ll}  
	2^{n_0-s} ( 2^{n_1} ( 2^{n_2} ( \cdots (2^{n_{i-1}} - \delta(i-1,i)) - \cdots ) - \delta(2, i) ) - \delta(1,i) ) & \mbox{if } i>0 \\
	1 & \mbox{if } i=0 \end{array} \right.  $$
where, as above, we take $l_{k+1}=-1$.  It is now straightforward to see:

\begin{theorem} \label{spine count}
For any $\tau$-function of length $N$, the number of truncated spines with this $\tau$-function is:
	$$\Spines(\tau) = \prod_{j=1}^N \SF(\tau, j).$$
\end{theorem}

\bigskip
\section{Complete algorithm}

We combine the results of the previous sections to produce an algorithm for the complete count of topological conjugacy classes for cubic polynomials with generic level $N$.

Fix a $\tau$-function of length $N$.  Recall that the markers of $\tau$ are the integers 
	$$\{m\in \{1, \ldots, N-1\}: \tau(m+1) < \tau(m)+1\},$$
and the marked levels are:
	$$\{l\geq 0: l = \tau^n(m) \mbox{ for some } n > 0 \mbox{ and marker } m\} \cup \{0\};$$
Let $L$ be the number of non-zero marked levels.  

For each $n\leq N$, the order of $n$ was defined in \S\ref{properties}; it satisfies $\tau^{\ord(n)}(n)=0$. For each marked level $l$, compute
	$$\mod(l) = \sum_{i=1}^l 2^{-\ord(i)}$$
and
	$$t(l) = \min\{n>0: n\mod(l) \in \N\}.$$
The quantity $\mod(l)$ represents the sum of relative moduli of annuli down to level $l$, while $t(l)$ is the numer of twists required to return that marked level to its original configuration.  We define
	$$T(\tau) = \max\{t(l): l \mbox{ is a marked level} \}$$
or set $T(\tau) = 1$ if $\tau$ has no marked levels.  The {\em twist factor} is defined by
	$$\TF(\tau) = \frac{2^L}{T(\tau)}.$$

From Theorem \ref{spine count}, the number of truncated spines with this $\tau$-function is:
	$$\Spines(\tau) = \prod_{j=1}^N \SF(\tau, j).$$
By \cite{DP:combinatorics}, the number of topological conjugacy classes associated to $\tau$ is then 
	$$\Top(\tau) = \Spines(\tau) \cdot \TF(\tau).$$
Combining these computations with Theorem \ref{tau extension}, it is straighforward to automate an inductive construction of all $\tau$-functions of length $N$, and we obtain an enumeration of all truncated spines and all topological conjugacy classes of generic level $N$.

\bigskip\bigskip
\def\cprime{$'$}

 \end{document}